\documentclass[11pt]{amsart}
\usepackage{amsfonts, amsmath, amssymb}
\usepackage{amscd}
\usepackage[pdftex, pdfstartview=FitH]{hyperref}
\usepackage[utf8]{inputenc}
\pdfoutput=1 

\input xy
\xyoption{all}

\newcommand{\R}{{\mathbb R}}

\newcommand{\Z}{{\mathbb Z}}
\newcommand{\RP}{{\mathbb {RP}}}
\newcommand{\CP}{{\mathbb {CP}}}
\newcommand{\HP}{{\mathbb {HP}}}
\newcommand{\vol}{{\rm vol}}

\newcommand{\geod}[1]{{\rm Geod}(#1)}

\newcommand{\length}{{\rm length}}
\newcommand{\slit}{{TM \! \setminus \! \mathbf{0}}}
\newcommand\restrict[1]{\raisebox{-.5ex}{$|$}_{#1}}

\newtheorem{theorem}{Theorem}[section]
\newtheorem{corollary}[theorem]{Corollary}
\newtheorem{lemma}[theorem]{Lemma}
\newtheorem{proposition}[theorem]{Proposition}

\newtheorem{innertheorem}{Theorem}
\newenvironment{maintheorem}[1]
    {\renewcommand\theinnertheorem{#1}\innertheorem}
    {\endinnertheorem}

\theoremstyle{definition}
\newtheorem{definition}[theorem]{Definition}
\newtheorem*{definition*}{Definition}

\title{Rigidity results for geodesically reversible Finsler metrics}

\author{J.C. \'Alvarez Paiva}

\address{J.C. \'Alvarez Paiva, U.M.R. CNRS 8524
         U.F.R. de Math\'ematiques, 59655 Villeneuve d'Ascq C\'edex, France}

\email{alvarez@math.univ-lille.fr}
\email{jalvarezpaiva@gmail.com}

\thanks{The author was partially supported by the grant ANR12-BS01-0009}



\begin{document}%
%%%%%%%%%%%%%%%%%

%%%%%%%%%%%%%%%%%%%%%%ABSTRACT%%%%%%%%%%%%%%%%%%%%%%%%%%%%%%%%%%%%%%%%%%
\begin{abstract} 
A Finsler metric is geodesically reversible if geodesics remain geodesics after a change of
orientation. Asymmetric norms on vector spaces and Funk metrics in the interior of convex bodies are examples of geodesically reversible metrics that are not necessarily sums of reversible metrics and closed $1$-forms. However,  there seem to be few such examples on closed manifolds. In this paper the theory of volumes and areas on Finsler spaces is applied to establish a number of rigidity theorems which partially explain this paucity of examples. These rigidity results settle some hitherto unsolved cases of Hilbert's fourth problem for asymmetric metrics.
\end{abstract}
%%%%%%%%%%%%%%%%%%%%%%%%%%%%%%%%%%%%%%%%%%%%%%%%%%%%%%%%%%%%%%%%%%%%%%%%
\maketitle
%%%%%%%%%%%%%%%%%%%%%%%%%%%%%%%%%%%
\begin{flushright}
    \begin{small}
        \parbox{4in} {{\it Slanders, sir: for the satirical slave says here that old men
                           have grey beards; that their faces are wrinkled; their eyes
                           purging thick amber and plum-tree gum; and that they have a
                           plentiful lack of wit, together with most weak hams: all which,
                           sir, though I most powerfully and potently believe, yet I hold it
                           not honesty to have it thus set down; for you yourself, sir,
                           should be old as I am, if, like a crab, you could go backward.}} \\
        --- Hamlet (Act II, Scene 2)
    \end{small}
\end{flushright}

%%%%%%%%%%%%%%%%%%%%%%%%%%%%%%%%%%%
%\tableofcontents

%%%%%%%%%%%%%%%%%%%%%%
\section{Introduction}
%%%%%%%%%%%%%%%%%%%%%%
The dual purpose of this paper is to study Hilbert's fourth problem in the compact setting as well as to explore the concept of geodesic reversibility introduced by Bryant in~\cite{Bryant:2006}.

We recall that Hilbert's fourth problem, in its succinct formulation due to Busemann~(\cite{Busemann:Hilbert}), asks {\it to construct and study all (possibly) asymmetric continuous metrics on open convex subsets of real projective $n$-space for which oriented projective lines are geodesics}. When the distances are symmetric, and the orientation of the lines is irrelevant, Busemann~(\cite{Busemann:1961b}) gave a simple integral-geometric construction for these {\it projective metrics}, while Pogorelov~\cite{Pogorelov:1979} and Szab\'o~\cite{Szabo:1986} showed that every symmetric projective metric can be obtained as a limit of those constructed by Busemann. The reader will find an enticing history of the problem in~\cite{Papadopoulos:2014}. 

In this paper we show that the problem of constructing all asymmetric projective metrics, which remains open in the case of open convex subsets of $\R^n$, has a surprisingly rigid answer on real projective $n$-space: 

\begin{maintheorem}{I} \label{spherical-space-forms}
A smooth projective Finsler metric on the real projective space $\RP^n$ is the sum of a reversible
projective Finsler metric and the differential of a smooth function on $\RP^n$. Moreover, 
 if $d$ is a continuous asymmetric projective metric on $\RP^n$, there exists a continuous function $f : \RP^n \rightarrow \R$ such that $d(x,y) - d(y,x) = f(y) - f(x)$ for every pair of points $x,y \in \RP^n$.
\end{maintheorem}

In other words, the {\it only} way to construct a continuous asymmetric 
distance function $d$ on $\RP^n$ for which all oriented projective lines are geodesics is to take a  symmetric distance function  $\Tilde{d}$  with the same property, find a continuous function 
$f: \RP^n \rightarrow \R$ such that 
$$
\Tilde{d}(x,y) > |f(y) - f(x)| 
$$
for all pairs of distinct points, and set
$$
d(x,y) = \Tilde{d}(x,y) + f(y) - f(x) .
$$

The two-dimensional case of the first part of Theorem~\ref{spherical-space-forms} is due to Gautier Berck who in 2012 privately communicated a proof---based on the range characterization of an integral-geometric transform---that a smooth projective Finsler metric on the $2$-sphere is the sum of a reversible Finsler metric and the differential of a function on $S^2$. Berck's proof is strictly two-dimensional, but armed with the knowledge that such rigidity was possible the author was able to find a different proof that works in all dimensions and applies to a much larger class of Finsler metrics.   

\begin{definition}[Bryant~\cite{Bryant:2006}]
A Finsler metric is said to be {\it geodesically reversible} if any geodesic can be reparametrized in an orientation-reversing way so as to remain a geodesic. 
\end{definition}

\begin{maintheorem}{II} \label{Zoll}
If $(M,F)$ is a smooth, geodesically-reversible Zoll Finsler manifold, then $F$ is the sum of a reversible Zoll Finsler metric and the differential of a smooth function on $M$.
\end{maintheorem}

A proof of this result was sketched in the preprint~\cite{Alvarez:2013} of which this paper is the final version. However, the underlying cause for this rigidity was still mysterious. Indeed, there are many examples of geodesically reversible or even projective Finsler metrics on open subsets of $\R^n$  (e.g., asymmetric norms in $\R^n$, the Funk metric in the interior of a convex body) that are not sums of reversible metrics plus closed $1$-forms, but on compact manifolds it seems that the only such examples are all flat. This gave rise to the project of understanding projective and other geodesically-reversible Finsler metrics on compact manifolds. Since the only manifolds with a flat projective structure (also known as a $(PGL(n+1,\R), \RP^n)$-structure) which are known to carry a $C^2$ projective Finsler metric are also those that carry a projective Riemannian metric (necessarily of constant curvature), this led to main problem tackled in~\cite{Alvarez-Barbosa:2018} and in the present work:  {\it to describe all $C^2$ projective Finsler metrics defined on compact or finite-volume space forms.}

Note that the Busemann-Pogorelov-Szab\'o solution of Hilbert's fourth problem in the reversible setting together with Theorem~\ref{spherical-space-forms} provides a solution for spherical space forms, even when  considering continuous distance functions. For compact Euclidean space forms the problem was solved by \'Alvarez Paiva and Barbosa Gomes in~\cite{Alvarez-Barbosa:2018} by showing that all such metrics are the sum of a flat Finsler metric and a closed $1$-form. In the hyperbolic case we have the following two results:

\begin{maintheorem}{III} \label{hyperbolic-space-forms-non-reversible}
 In two and three dimensions, every $C^2$ projective Finsler metric on a hyperbolic space form of finite volume is the sum of a multiple of the hyperbolic metric and a closed $1$-form.   
\end{maintheorem}

 \begin{maintheorem}{IV} \label{hyperbolic-space-forms-reversible}
If $F$ is a reversible $C^2$ projective Finsler metric on a hyperbolic space form of finite volume, then $F$ is a multiple of the hyperbolic metric.
\end{maintheorem}

Like in the case of Euclidean space forms treated in~\cite{Alvarez-Barbosa:2018}, the regularity of the metric is important: examples by Busemann, Salzmann, and 
Bonahon~(\cite{Busemann-Salzmann:1965} and~\cite{Bonahon:1993}) show that there are plenty of non-proportional, continuous, symmetric projective metrics on hyperbolic surfaces. However, like in the case of spherical space forms, the results for hyperbolic space forms of finite volume follow from more general results on geodesic-reversibility and projective rigidity that depend on the dynamics of the geodesic flow.  

\begin{maintheorem}{V}\label{geodesic-reversibility}
Let $M$ be a two- or three-dimensional manifold and let $F$ be a $C^2$ geodesically reversible Finsler metric on $M$. If the geodesic flow does not admit any non-trivial continuous integral of motion, then the metric $F$ is necessarily the sum of a reversible metric and a closed 1-form.
\end{maintheorem}

\begin{maintheorem}{VI} \label{transitive-geodesic}
Let $F_1$ and $F_2$ be two reversible $C^2$ Finsler metrics defined on a connected manifold $M$ and which have the same unparametrized geodesics. If the geodesic flow of either metric does not admit any non-trivial continuous integral of motion, then the metric $F_2$ is necessarily a  constant multiple of $F_1$.
\end{maintheorem}

Some examples of Finsler manifolds whose geodesic flow does not admit any non-trivial continuous integrals of motion are (1) compact Finsler manifolds with negative flag curvature (Foulon~\cite{Foulon:1992}), (2) compact Finsler surfaces of genus greater than one without conjugate points (Barbosa Gomes and Ruggiero~\cite{Barbosa-Ruggiero:2013}), and (3) compact Finsler manifolds without conjugate points and whose universal cover satisfies the uniform Finsler visibility condition (Chimenton, Barbosa Gomes, and Ruggiero~\cite{Barbosa-Chimenton-Ruggiero:2020}). 

Theorem~\ref{transitive-geodesic} is the generalization to reversible Finsler metrics of a result by Matveev and Topalov (see Corollary~1 in~\cite{Matveev-Topalov:2000, Matveev-Topalov:2003}). Note however that both the local and global study of Riemannian metrics with the same unparametrized geodesics is full of rigidity phenomena that do not carry over to Finsler metrics.   For instance, a simple consequence of Corollary~2 in~\cite{Matveev-Topalov:2003} is that {\it two Riemannian metrics in a connected manifold that have the same unparametrized geodesics and agree on some open set are necessarily identical.} This is false for Finsler metrics, even reversible ones: the Busemann construction of projective metrics on the sphere allows us to construct an abundance of non-isometric examples of smooth reversible Finsler metrics all of whose geodesics are great circles and which agree with the standard Riemannian metric in the complement of any fixed neighborhood of  the equator (see Section~\ref{Zoll-metrics}). Nevertheless, we do have the following rigidity phenomenon:

\begin{maintheorem}{VII} \label{continuation-principle}
Let $F_1$ and $F_2$ be two $C^2$ Finsler metrics with the same unparametrized geodesics. If the velocity field of every geodesic passes through some open set in the tangent bundle where $F_1$ equals $F_2$, then the metrics agree everywhere modulo the addition of a closed $1$-form. 
\end{maintheorem}

An interesting and surprising aspect of the proofs of Theorems~\ref{Zoll}, ~\ref{geodesic-reversibility}, and~\ref{transitive-geodesic}  is that they are 
based on the study of volumes and areas on Finsler spaces (see~\cite{Thompson:1996} for a thorough account of the convex-geometric origins of the subject and~\cite{Alvarez-Thompson:2004} for a more Finsler-oriented and updated account). In the final analysis, Theorem~\ref{Zoll} rests on the Brunn-Minkowski theorem, Theorem~\ref{geodesic-reversibility} on Chakerian's characterization of bodies of constant width and constant brightness in relative geometry~(\cite{Chakerian:1967}), and Theorem~\ref{transitive-geodesic} on the solution of the Minkowski problem. Despite---or perhaps due to---the interplay between symplectic, integral and convex geometry on which they are based, the proofs are exceedingly simple. In Sections~\ref{volumes} and~\ref{Hilbert-form} we review some basic concepts and develop an adequate language in which to describe this interplay.
Section~\ref{Hilbert-form} also contains the proof of Theorem~\ref{continuation-principle}. Sections~\ref{Zoll-metrics} and~\ref{geodesic-rigidity} contain the proofs of Theorems~\ref{spherical-space-forms} to~\ref{transitive-geodesic}.  

\medskip
\noindent {\bf Acknowledgments.} The author thanks G.~Berck, A.~Thompson, J.~Barbosa~Gomes, S.~Tabachnikov, B.~McKay, and V.~Matveev for the many interesting discussions that have in some form or another found their way to the following pages. I'm specially grateful to R.~Howard and D.~Ryabogin for their timely answers to my naive questions about bodies of constant width and constant brightness. The author is also indebted to R.~Ruggiero and R.~Bryant for their comments on the first version of this work.

%%%%%%%%%%%%%%%%%%%%%%%%%%%%%%%%%%%%%%%%%%%%%%%%%%%%%%%%%%%%%
\section{Volumes and areas in Finsler spaces} \label{volumes}
%%%%%%%%%%%%%%%%%%%%%%%%%%%%%%%%%%%%%%%%%%%%%%%%%%%%%%%%%%%%%

A {\sl continuous Finsler metric} on a manifold $M$ is a continuous function $F:TM \rightarrow [0,\infty)$ such that its restriction to each tangent space $T_x M$ $(x \in M)$ is a (possibly) asymmetric norm. In order to study things like geodesic flows, curvature, and differential invariants, it is usual to ask for more regularity and a stricter notion of convexity:

\begin{definition}
Let $V$ be a finite-dimensional vector space over the reals. An asymmetric norm $\|\cdot\| : V \rightarrow [0,\infty)$ is said to be a {\sl Minkowski norm} if  outside the origin the function $\|\cdot\|^2$ is $C^2$ and its Hessian is positive-definite. A continuous function 
$F : T M \rightarrow [0, \infty)$ defined on the tangent bundle of a smooth manifold $M$ will be said to be a $C^k$ $(k \geq 2)$ {\sl Finsler metric} if it is $C^k$ outside the zero section and its restriction to each tangent space is a Minkoswki norm.
\end{definition}

Note that the restriction of the Finsler metric $F$ to a tangent space $T_xM$ is the support function of the {\sl unit co-disc} 
$$
D_x^* M := \{\xi_x \in T_x^* M: \xi_x(v_x) \leq 1 \text{ whenever } F(v_x) \leq 1 \}.
$$
The {\sl unit co-disc bundle} $D^* M \subset T^*M$ of a Finsler metric on a manifold $M$ is the union of all the co-discs $D_x^* M$, $x \in M$, seen as a disc bundle over $M$.  

\begin{definition}
The {\sl Holmes-Thompson volume} of a Finsler metric on a manifold $M$ of dimension $n$ is the symplectic volume of its unit co-disc bundle divided by the volume of the $n$-dimensional Euclidean unit disc. The {\sl area} of a submanifold is defined as the Holmes-Thompson volume of the submanifold provided with its inherited metric.
\end{definition}

It is also possible to define the Holmes-Thompson volume by considering its volume density---the function $\Phi : \bigwedge^n T M  \to [0,\infty)$ that measures the volume of parallelotopes formed by vectors tangent to $M$---:  if $v_1, \ldots, v_n $  is a basis of $T_xM$, the $n$-vector 
$v_1 \wedge \cdots \wedge v_n$ is a volume form on $T_x^* M$. Denoting the volume of the $n$-dimensional Euclidean unit disc by $\varepsilon_n$, we define
$$
\Phi(v_1 \wedge \cdots \wedge v_n) :=  \frac{1}{\varepsilon_n} \int_{D_x^* M} |v_1 \wedge \cdots \wedge v_n| .
$$
The equivalence between both definitions follows immediately from the computation of the 
pushfoward of the symplectic volume form on $D^*M$ onto $M$ under the canonical projection 
$\pi : D^* M \rightarrow M$.

From this second description it is clear that adding a $1$-form to a Finsler metric does not alter its volume density: adding the form merely translates the unit co-discs in each cotangent space. A more interesting property of the Holmes-Thompson volume is its behavior with respect to {\sl central symmetrization}: 

\begin{definition}
The central symmetrization of a Finsler metric  $F$ is the Finsler metric
defined by $F^S(v_x) = (F(v_x) + F(-v_x))/2$. 
\end{definition}

\begin{proposition}\label{central-symmetrization}
The Holmes-Thompson volume density of a Finsler metric on a manifold $M$ is less than or equal to that of its central symmetrization. Moreover, the volume densities are equal if and only if $F = F^S + \beta$ for some $1$-form $\beta$ on $M$.   
\end{proposition}

\proof 
Note that if $D_x^* M \subset T^*M$ is the convex body supported by the restriction of $F$ to $T_xM$, the restriction of $F^S$ to this vector space is the support function of the symmetrized body 
$\Delta D_x^*M := (D_x^* M - D_x^* M)/2$. A well-known consequence of the Brunn-Minkowski inequality is that the volume of a convex body $K$ is less than or equal to that of its symmetrization $\Delta K = (K - K)/2$ with equality if and only if $K$ is symmetric with respect to a (unique) point $p$ in its  interior. Applying the theorem at every cotangent space $T^*_x M$, $x \in M$, concludes the proof. Note that the $1$-form will be smooth if the bodies $D_x^*M$ depend smoothly on the base point.
\qed

The idea in the preceding proof yields a somewhat more general result:
 
\begin{proposition}\label{Brunn-Minkowski}
Let $F_1$ and $F_2$ be two continuous Finsler metrics on a manifold $M$. If the Holmes-Thompson densities of $(M,F_1)$ and $(M,F_2)$ are identical, then they are both
less than or equal to the Holmes-Thompson volume density of the average metric $(F_1 + F_2)/2$. Moreover, equality holds if and only if $F_2 = F_1 + \beta$  for some $1$-form $\beta$ on $M$.  
\end{proposition}

\proof[Sketch of the proof.] Notice that at each point of $M$ the unit codisc of the average metric is the Minkowski sum of the the codiscs of the metrics $F_1/2$ and $F_2/2$ at that point. Applying the Brunn-Minkowski inequality and its equality case concludes the proof.
\qed

%%%%%%%%%%%%%%%%%%%%%%%%%%%%%%%%%%%%
\subsection{Areas in Finsler spaces}
%%%%%%%%%%%%%%%%%%%%%%%%%%%%%%%%%%%%

The extrinsic description of areas in Finsler spaces requires us to define the area of each $k$-dimensional parallelotope in $T_x M$ at every point $x \in M$. Given a decomposable $k$-vector $a$ formed by taking the exterior product of linearly independent tangent vectors $v_1,\ldots,v_k$
in $T_xM$, consider the subspace $\langle a \rangle$ spanned by these vectors together the natural inclusion $i_a : \langle a \rangle \rightarrow T_xM$.  The dual map  $i_a^* : T_x^* M \rightarrow \langle a \rangle^*$ projects the unit co-disc $D_x^*M \subset T_x^*M$ 
to a convex body $i_a^*(D_x^*M)$ in $\langle a \rangle^*$. Since $a$ is a volume form in  $\langle a \rangle^*$, we are able to define the {\sl Holmes-Thompson $k$-area density}
$$
\varphi_k(a_x) = \frac{1}{\varepsilon_k} \int_{i_a^*(D_x^*M)} |a| .
$$

When the Finsler manifold is $\R^n$ provided with a norm with unit ball $D$ the Holmes-Thompson $k$-area density can be identified with the {\it $k$-th projection function} of the polar body $D^\circ$~(\cite[Chapter VI]{Thompson:1996}). Since this provides some intuition and allows us to connect the theory of areas in Finsler spaces with important results in convex geometry, we will go carefully over the identification. 

Consider $\R^n$ provided with its standard inner product together with a convex body $D$ containing the origin in its interior. The dual space $\R^{n*}$ is identified with $\R^n$ via the inner product, and the dual body $D^*$ is identified with the polar body
$$
D^\circ := \{u \in \R^n : u \cdot v \leq 1  \text{ whenever } v \in D \} .
$$
Given orthonormal vectors $e_1,\ldots,e_k$ in $\R^n$, we set $a := e_1 \wedge \cdots \wedge e_k$ 
and denote the subspace they span by $\langle a \rangle$. Using the restriction of the inner product
on $\R^n$ to $\langle a \rangle$ we can identify this last space with its dual space 
$\langle a \rangle^*$. With these identifications, the map  $i_a^* : \R^{n*} \rightarrow \langle a \rangle^*$ used
in the definition of the $k$-area density becomes the orthogonal projection from $\R^n$ to  $\langle a \rangle$. {\it In conclusion, taking all these identifications into account, $\varphi_k(a)$ is simply the Euclidean area of the orthogonal projection of $D^\circ$ onto  $\langle a \rangle$ divided by the area of the $k$-dimensional unit disc in Euclidean space.}  

This geometric intuition makes it obvious that adding a $1$-form to $F$ or considering the {\it reverse metric} $F_{-}$, defined by $F_{-}(v_x) = F(-v_x)$, does not change the $k$-area densities. The following results also follow from the relation between the Holmes-Thompson $k$-area densities and $k$-th projection functions in convex geometry. An excellent reference for what follows is Chapter~3 of~\cite{Gardner:1995}. 

\begin{proposition}\label{areas-determine-metric}
Given a Finsler metric $F$ on a manifold $M$, there exists a unique reversible Finsler metric $F^B$ on $M$, the areal symmetrization of $F$, that defines the same $(n-1)$-area density. 
\end{proposition}

\proof
A classic result in convex geometry states that given a convex body $K \subset \R^n$, there is a unique $0$-symmetric body $\nabla K$ called the {\it Blaschke body} of $K$ such that the areas of orthogonal projections of $K$ and $\nabla K$ onto any fixed hyperplane are the same. Applying this 
construction to every unit co-disc we obtain a  reversible Finsler metric $F^B$ whose restriction to every tangent space $T_xM$ is the support function of the Blaschke body of $D_x^*M$. By construction $F$ and $F^B$ have the same $(n-1)$-area density.
\qed

\begin{proposition}
The Holmes-Thompson volume density of a Finsler metric $F$ on a manifold $M$ is less than
or equal to that of its areal symmetrization $F^B$. Moreover, the volume densities are equal if 
and only if $F = F^B + \beta$ for some $1$-form $\beta$ on $M$. 
\end{proposition}

\proof
By the Kneser-S\"uss inequality (Theorem~3.3.9 in~\cite[p. 109]{Gardner:1995}), if $K \subset \R^n$ is a convex body, its volume is less than or equal to the volume of its Blaschke body $\nabla K$ and equality holds if and only if $K$ is centrally symmetric. The result follows by applying this inequality to every unit co-disc of the Finsler manifold $(M,F)$ in the same way we applied the Brunn-Minkowski inequality in Proposition~\ref{central-symmetrization}.
\qed

\begin{proposition}\label{constant-width-constant-brightness}
Given a $C^2$ Finsler metric $F$ on a three-dimensional manifold $M$, its central symmetrization  $F^S$ is conformal to its areal symmetrization $F^B$ (i.e., $F^S = e^f F^B$, where $f$ is a continuous function on $M$) if and only if $F = F^S + \beta$ for some $1$-form $\beta$ on $M$.
\end{proposition}

\proof
Let us fix an arbitrary point $x \in M$ and consider what the hypotheses tell us about the convex body $D_x^* M$. Since $F^S$ and $F^B$ are assumed to be conformal, their restrictions to $T_xM$ are multiples and
thus the central symmetrization and the Blaschke body of $D_x^* M$ are dilates: 
$\Delta D_x^*M = \lambda \nabla D_x^*M$ for some $\lambda > 0$. 

Note that, by definition of the symmetrized body, $D_x^* M$ is a body of constant width with respect to $\Delta D_x^*M$. Likewise, by definition of the Blaschke body, $D_x^* M$ is a body of constant brightness with respect to $ \nabla D_x^*M$. Since $\Delta D_x^*M$ and $\nabla D_x^*M$ are dilates, this implies that $D_x^*M$ is a body of constant width and constant brightness relative to the body $\Delta D_x^* M$.  A well-known result of
of Chakerian~\cite{Chakerian:1967} (see also~\cite{Howard:2006}) states that a convex body that is of constant width and constant brightness relative to a centrally symmetric {\it gauge body} whose boundary is  $C^2$ and quadratically convex is necessarily obtained from the gauge body by translation and dilation. Considering that $F^S$ is a $C^2$ Finsler metric, $\Delta D_x^*M$ is precisely such a gauge body. Chakerian's result then implies that $D_x^*M$ is centrally symmetric for every point $x \in M$, and hence $F$ is the sum of a reversible metric and a $1$-form.
\qed

The generalization of Chakerian's result in~\cite{Chakerian:1967} to dimensions greater than three is a long-standing open problem in convex geometry. As we shall see in Section~\ref{geodesic-rigidity}, should the generalization prove true, Theorem~\ref{geodesic-reversibility} would also generalize to arbitrary dimension.

%%%%%%%%%%%%%%%%%%%%%%%%%%%%%%%%%%%%%%%%%%%%%%
\subsection{Basic integral-geometric formulas}
%%%%%%%%%%%%%%%%%%%%%%%%%%%%%%%%%%%%%%%%%%%%%%

Let $(M,F)$ be a $C^2$ Finsler manifold and let $S^* M := \partial D^*M$ be its {\it unit co-sphere bundle}. The pullback of the canonical $1$-form $\alpha$ in $T^*M$ to the co-sphere bundle is a $C^1$ contact form and its Reeb vector field defines the geodesic (local) flow of $(M,F)$. 

Assuming the metric is $C^3$ and the space of oriented geodesics of $(M,F)$ is itself a manifold $\geod{M}$,  it is well-known (cf.~\cite[p.~48--49]{Arnold-Givental:1990}) that $\geod{M}$  carries a symplectic $2$-form $\omega_F$ defined as follows: if  
$$
\begin{CD}
S^*M  @>i>>  T^*M \\
@V{\pi}VV          \\
\geod{M}
\end{CD}
$$
are the canonical projection onto $\geod{M}$ and the canonical inclusion
into $T^{*}M$, $\omega_F$ is the unique $2$-form satisfying the equation $\pi^{*}\omega_F = i^{*}d\alpha$. This is standard symplectic reduction and, as it is explained in Theorem~5.3.23 in~\cite[p.~416]{Abraham-Marsden}, it requires the symplectic form $d\alpha$ to 
be $C^1$, which in turn requires the metric $F$ to be $C^3$. However, with very little effort we can get by with less. Employing the notation above, we have the following

%%%%%%%%%%%%%%%%%%%%%%%%%% I'm here %%%%%%%%%%%%%%%%%%%%%%%%%%%%%%%%%%%%%%%%%%%%%%%%%%%%%%%%%%%%

\begin{theorem}\label{space-of-geodesics}
Let $(M,F)$ be a $C^2$ Finsler manifold and assume the space of geodesics of $(M,F)$ is itself
a manifold $\geod{M}$. There is a unique continuous, maximally non-degenerate $2$-form $\omega_F$
on $\geod{M}$ that satisfies the equation $\pi^{*}\omega_F = i^{*}d\alpha$.
\end{theorem}

The key point is to show that the continuous form $i^{*}d\alpha$ is invariant under the geodesic 
flow without being able to take Lie derivatives. Since a slight generalization of this is also the key to Theorem~\ref{continuation-principle}, we encapsulate what we need in one result. 

\begin{lemma}\label{surfaces-of-section}
Let $P$ be a contact manifold with $C^1$ contact form $\alpha$ and let  $S_1, S_2 \subset P$ be two hypersurfaces that are transverse to its Reeb vector field. If every Reeb orbit that intersects $S_1$ also intersects $S_2$ at a later time, then the map  $\tau: S_1 \rightarrow S_2$  that sends each point $p \in S_1$ to the first crossing of $S_2$ and the
Reeb orbit with initial condition $p$ is an isomorphism between $(S_1, d\alpha)$ and $(S_2,d\alpha)$. \end{lemma}

\proof
It is enough to prove that the integral of $d\alpha$ over every sufficiently small $2$-disc $R \subset S_1$ equals the integral of $d\alpha$ over its image $\tau(R)$. In order to do this, consider
the cylinder $C \subset P$ formed by the segments of Reeb orbits joining the points in the boundary of $R$ with those in the boundary of $\tau(R)$. Taking into account that tangent vectors to  Reeb orbits are in the kernel of $d\alpha$ and applying Stokes formula (twice) we have that 
$$
0 = \int_C d\alpha = \int_{\partial R} \alpha - \int_{\partial \tau(R)} \alpha =
\int_R d\alpha - \int_{\tau(R)} d\alpha .
$$
\qed

\proof[Proof of Theorem~\ref{space-of-geodesics}]
Given a sufficiently small open set $U \subset \geod{M}$ take a hypersurface $\Sigma \subset S^*M$ that is transverse to the flow and projects diffeomorphically onto $U$. Define $\omega_F$ by 
setting $\pi^* \omega_F$ to be equal to $d\alpha$ on $\Sigma$. By Lemma~\ref{surfaces-of-section}, the choice of transverse hypersurface is immaterial. The form $\omega_F$ is maximally non-degenerate because  $\pi^* \omega_F^{n-1}$ is equal to the nowhere-zero form $d\alpha^{n-1}$ on $\Sigma$. 
\qed

Among the many examples of Finsler manifolds whose space of geodesics is
a manifold we find convex neighborhoods in Finsler spaces,
Hadamard Riemannian manifolds (for both of these examples see
\cite{Ferrand:1997}), Zoll metrics (\cite{Besse:1978}), Minkowski geometries, Hilbert geometries,
Funk geometries and, more generally, projective Finsler metrics (\cite{Alvarez:2005}
and \cite{Alvarez-Fernandes:1998}). For these manifolds we have the following
integral-geometric formulas:

\begin{theorem}[\cite{Alvarez-Berck:2006}]\label{crofton-hypersurface}
Let $M$ be an $n$-dimensional $C^2$ Finsler manifold with manifold of
geodesics $\geod{M}$. If $\omega_F$ is the induced $2$-form on $\geod{M}$
and  $N \subset M$ is any immersed hypersurface, then
$$
\vol_{n-1}(N) = \frac{1}{2\epsilon_{n-1}(n-1)!} \, \int_{\gamma \in \geod{M}}
\#(\gamma \cap N) |\omega_F^{n-1}| ,
$$
where $\epsilon_{n-1}$ is the volume of the Euclidean unit ball of dimension
$n-1$.
\end{theorem}

\begin{proposition}[Santal\'o formula]
Let $F$ be a $C^2$ Finsler metric on an $n$-dimensional manifold $M$. If $f : S^*M \rightarrow \R$ is a continuous function on its unit co-sphere bundle, then
$$
\int_{S^*M} f \alpha \wedge d\alpha^{n-1} = 
\int_{l \in \geod{M}} \! \! \left( \int_l f \alpha \right) \omega_F^{n-1} .
$$
In particular, taking $f \equiv 1$, 
$$
\vol_n(M,F) = \frac{1}{\varepsilon_n n!} \int_{S^*M} \! \! \alpha \wedge d\alpha^{n-1} 
= \frac{1}{\varepsilon_n n!} \int_{l \in \geod{M}} \! \! \length(l) \, \omega_F^{n-1} .
$$
\end{proposition}

\proof
Using that the tangents to the fibers of the projection $\pi : S^*M \rightarrow \geod{M}$ are in
the kernel of the form $d\alpha^{n-1} = \pi^* \omega_F^{n-1}$ we obtain the identity
$$
\pi_*  \left( f \alpha \wedge d\alpha^{n-1} \right) = \left( \int_l f \alpha \right) \omega_F^{n-1} ,
$$
which immediately implies the result. 
\qed

%%%%%%%%%%%%%%%%%%%%%%%%%%%%%%%%%%%%%%%%%%%%%%%%%%%%%%%%%%%%%%%%%%%%%%%%%%
\section{The Hilbert form and geodesic reversibility} \label{Hilbert-form}
%%%%%%%%%%%%%%%%%%%%%%%%%%%%%%%%%%%%%%%%%%%%%%%%%%%%%%%%%%%%%%%%%%%%%%%%%%

A common feature of differential $1$-forms and Finsler metrics is that the result of integrating them over a continuously differentiable curve is independent of the parameterization as long as the orientation of the curve is unchanged. In the case of Riemannian and reversible Finsler metrics, the integral is also invariant under a change of orientation.  General integrands that satisfy these properties are known, respectively, as $1$-densities and even $1$-densities:

\begin{definition}
A continuous function $L : T M \rightarrow \R$ defined on the tangent bundle of a smooth manifold $M$ will be said to be a $1$-{\sl density} of class $C^k$ $(k \geq 1)$
if it is  $C^k$ outside the zero section and positively homogeneous of degree $1$ in the velocities (i.e., $L(tv_x) = tL(v_x)$ for $t > 0$). A $1$-density $L$
will be called  {\sl even} if it  satisfies $L(-v_x) = L(v_x)$ and  {\sl odd} if it satisfies $L(-v_x) = -L(v_x)$.
\end{definition}

Note that $C^k$ Finsler metrics form a convex cone inside the vector space of $C^k$ densities.
Moreover, every Finsler metric $F$ can be uniquely decomposed into odd and even densities:  if $F_{-}$ denotes the reverse of the Finsler metric $F$, we can write
$$
F = \frac{1}{2}(F - F_{-})  + \frac{1}{2}(F + F_{-}) =:  F^{\mathcal{O}} + F^S  .
$$ 

%%%%%%%%%%%%%%%%%%%%%%%%%%%%%%%
\subsection{The Hilbert Form}
%%%%%%%%%%%%%%%%%%%%%%%%%%%%%%%
In what follows it will be convenient to adopt a manifestly covariant formalism for the calculus of variations due to Hilbert (see the text of his 23rd problem in~\cite{Hilbert:problems}). In modern terminology it is quite simple (cf.~\cite{Chern:1996}): given a $C^2$ $1$-density $L$ on $M$ we  define a bundle map $\delta L : \slit \rightarrow T^* M$ from the tangent bundle of $M$ minus its zero section to the cotangent bundle of $M$ 
by the formula
$$
\delta L(v_x) \cdot w_x := \frac{d}{dt}L(v_x + tw_x)\restrict{t = 0} \, .
$$
The map $\delta L$ is used to pull back the canonical $1$-form $\alpha$ on the cotangent bundle to a 
$1$-form $\alpha_{{}_L} := (\delta L)^* \alpha$ on $\slit $.  This is {\it the Hilbert form} associated to $L$.  In local coordinates
$$
\alpha_{{}_L} = \frac{\partial L}{\partial v_1} \, dx_1 + 
\cdots + \frac{\partial L}{\partial v_n} \, dx_n .
$$

Two immediate consequences of this local expression are that (1) the Hilbert $1$-form $\alpha_L$ is homogeneous of degree zero and so descends naturally onto the space of tangent directions $STM := (\slit) / \R^+$, and (2) that  the Hilbert form of a differential $1$-form $\beta \in \Omega^1(M)$, seen as the odd $1$-density $(x,v) \mapsto \beta_x(v)$, is just the pullback of  $\beta$ to $\slit$ under the canonical projection $\pi : \slit  \rightarrow M$.

Using Euler's formula for homogeneous functions we see at once that $L(x,v) = \alpha_{(x,v)}(v)$. In order to be precise, we need to interpret the symbol $v$ in two different ways: as a vector in $T_x M$ and as a vector in $T_{(x,v)}(\slit)$ (i.e., as the velocity  vector at $t = 0$ of the curve
$(x,v + tv)$ in the tangent bundle of $M$). It follows that for any continuously differentiable curve $\gamma$ with velocity curve $\gamma'$
$$
 \int_\gamma L = \int_{{\gamma}'} \alpha_{{}_L} .
$$

Using the Hilbert form, the Euler-Lagrange equations  for $1$-densities can be elegantly formulated as follows: 

\begin{theorem}\label{first-variation}
A smooth curve $\gamma : \R \rightarrow M$ with nowhere vanishing velocity is extremal for the variational problem $\gamma \mapsto \int_\gamma L$ 
if and only if its second derivative $\gamma''$, seen as  the velocity of the curve $\gamma'$, satisfies the equation
$$
d \alpha_{{}_L} (\gamma''(t), \cdot) = 0
$$
for all values of $t$. 
\end{theorem}

\proof
Given  a smooth, compactly-supported variation $\gamma_s$ of the curve $\gamma = \gamma_0$ we wish to compute the derivative at $s = 0$ of the 
function 
$$
s \longmapsto \int_{\gamma_s} L = \int_{\gamma'_s} \alpha_{{}_L}  .
$$
Writing $\gamma'_s = \phi_s(\gamma')$, where $\phi_s : \slit \rightarrow \slit$ is a smooth isotopy that equals the identity outside a compact set, we 
have that 
$$
\frac{d}{ds}\left( \int_{\gamma'_s} \alpha_{{}_L} \right) \restrict{s = 0} = \frac{d}{ds} \left( \int_{\phi_s(\gamma')} \alpha_{{}_L} \right) \restrict{s = 0} = 
\int_{\gamma'} \mathcal{L}_X \alpha_{{}_L} .
$$
Here $\mathcal{L}_X$ is the Lie derivative with respect to the  vector field on $\slit$ defined by $X(v_x) = d\phi_s(v_x)/ds\restrict{s = 0}$. Using Cartan's formula for the Lie derivative and assuming that either $\gamma$ is periodic or that the vector field $X$ is supported away from the endpoints of $\gamma'$,
$$
\int_{\gamma'} \mathcal{L}_X \alpha_{{}_L} = \int_{\gamma'} d\alpha_{{}_L}(X, \cdot) + d(\alpha_{{}_L}(X)) =  \int_{\gamma'} d\alpha_{{}_L}(X, \cdot) \, .
$$
If $\gamma$ is an extremal, this quantity must be zero for an arbitrary compactly-supported vector field $X$ and, therefore, we have that 
$d \alpha_{{}_L} (\gamma''(t), \cdot) = 0$ for all values of $t$. 
\qed

A direct consequence of this theorem---and arguably the most important remark in the study 
of inverse problems for $1$-densities and Finsler metrics---is the following {\it principle of superposition:}

\begin{corollary}
Given a family of smooth oriented curves on a manifold $M$, the set of $1$-densities on $M$ for which all the curves in the family are extremal is closed under linear combinations. In particular, the set of Finsler metrics on $M$ for which all the curves in the family are geodesic is closed under positive linear combinations.
\end{corollary}

\proof
It is enough to notice that the Hilbert form depends linearly on the $1$-density:
if $L_1$ and $L_2$ are two $1$-densities, then for any real numbers $\lambda$ and $\mu$, 
$\alpha_{{}_{\lambda L_1 + \mu L_2}} = \lambda \alpha_{{}_{L_1 }} + \mu \alpha_{{}_{L_2}}$. If a smooth curve $\gamma : \R \to M$ is an extremal 
for both $L_1$ and $L_2$, then 
$$
d\alpha_{{}_{\lambda L_1 + \mu L_2}} (\gamma'', \cdot) =  \lambda d\alpha_{{}_{L_1 }} (\gamma'', \cdot)  + \mu  d\alpha_{{}_{L_2}}  (\gamma'', \cdot)  = 0 ,
$$
which means it is also an extremal for $\lambda L_1 + \mu L_2$.
\qed

\begin{corollary}
If $F$ is a geodesically reversible $C^2$ Finsler metric, then $F$, its reverse metric $F_{-}$, and its central symmetrization $F^S$ all have the same unparametrized geodesics. Moreover, the geodesics
of $F$ are also extremals of its odd part $F^\mathcal{O}$.  
\end{corollary}

Theorem~\ref{first-variation} also implies that every smooth curve with never-vanishing velocity is an extremal for the $1$-density $L$ if and only if $d \alpha_{{}_L}$ 
vanishes identically. This happens if and only if $L$ is a closed differential form (cf.~\cite[p. 51--59]{Giaquinta-Hildebrandt:vol1}):

\begin{proposition}\label{null-lagrangian}
The exterior differential of the Hilbert form of a  $1$-density $L : T M \rightarrow \R$ vanishes identically if and only if $L$  is a closed $1$-form on $M$.
\end{proposition}

\proof
This is simplest to see in local coordinates:
$$
d \alpha_{{}_L} = d \left( \sum_{i=1}^n \frac{\partial L}{ \partial v_i} \, dx_i \right) = \sum_{1 \leq i,j \leq n} \frac{\partial^2 L}{ \partial v_j \partial v_i} \, dv_j \wedge dx_i
+ \sum_{1 \leq i,j \leq n} \frac{\partial^2 L}{\partial x_j \partial v_i} \, dx_j \wedge dx_i .
$$
If $d \alpha_{{}_L}$ vanishes identically, then both
$$
 \sum_{1 \leq i,j \leq n} \frac{\partial^2 L}{\partial v_j \partial v_i} \, dv_j \wedge dx_i    \  \  \ \hbox{\rm and }   \ 
  \sum_{1 \leq i,j \leq n} \frac{\partial^2 L}{\partial x_j \partial v_i} \, dx_j \wedge dx_i 
$$
vanish identically.  The vanishing of the first expression implies that $L$ is a $1$-form (i.e., that $L$ is linear in the velocities), the  vanishing of
the second  says that this form is closed. 
\qed

This in turn leads to a simple proof of the following key result (cf. Crampin~\cite{Crampin:2005}). 

\begin{corollary}\label{adding-form}
On a manifold $M$ consider a Finsler metric $F$ and a $1$-form $\beta$. The extremals of $F + \beta$ coincide with those of $F$ if and only if $\beta$ is closed. Moreover, if $F$ is reversible, then $F + \beta$ is geodesically reversible if and only if $\beta$ is closed.
\end{corollary}

\proof
Notice that $d\alpha_{{}_{F + \beta}} = d\alpha_{{}_{F}} + \pi^* d\beta$, where $\pi : \slit  \rightarrow M$ is the canonical projection.  It is clear
from Theorem~2.3 that if $\beta$ is closed, $F$ and $F + \beta$ have the same extremals. In order to prove the converse, denote the geodesic 
vector field of $(M,F)$ by $X$ and remark that, by hypothesis, it is in the kernel of  both $d\alpha_{{}_{F}}$ and $d\alpha_{{}_{F + \beta}}$. It follows 
that $(\pi^* d\beta) (X(v_x),\cdot) = d\beta(v_x, \cdot)$ is identically zero for every vector $v_x \in \slit$ and, therefore, $d\beta$ is identically zero.

To establish the second part of the corollary recall that a Finsler metric $L$ is geodesically reversible if and only if $L$ and its reverse have the same extremals. Since $F$ is reversible and $\beta$ is a $1$-form, it follows that $F + \beta$ is geodesically reversible if and only if $F + \beta$ and $F - \beta$ have the same extremals. This implies that the Finsler metric $F$ and the $1$-density $F + \beta$ have the same extremals and $\beta$ must be a closed $1$-form.
\qed

%%%%%%%%%%%%%%%%%%%%%%%%%%%%%%%%%%%%%%%%%%%%
\subsection{Basic geodesic rigidity results}
%%%%%%%%%%%%%%%%%%%%%%%%%%%%%%%%%%%%%%%%%%%%
The two rigidity results that follow are basic in the sense that they depend only on very simple and general assumptions on maximally-degenerate, closed $2$-forms on odd dimensional manifolds. Unlike the other results in this paper, the fact that these closed $2$-forms are the exterior differentials
of Hilbert $1$-forms of Finsler metrics is not fundamental.

\begin{maintheorem}{VII} 
Let $F_1$ and $F_2$ be two $C^2$ Finsler metrics with the same unparametrized geodesics. If the velocity field of every geodesic passes through some open set in the tangent bundle where $F_1$ equals $F_2$, then the metrics agree everywhere modulo the addition of a closed $1$-form. 
\end{maintheorem}

\proof
Consider the manifold of tangent directions $STM$ together with the Hilbert forms $\alpha_{F_1}$
and $\alpha_{F_2}$, which we consider as $C^1$ contact forms on this space. In light of Proposition~\ref{null-lagrangian}, if we prove that $d\alpha_{F_1}$ and $d\alpha_{F_2}$ are identical, we will have proved that $F_1$ and $F_2$ differ by a closed $1$-form. 

Note that since the metrics $F_1$ and $F_2$ have the same unparametrized geodesics, the kernels
of the forms  $d\alpha_{F_1}$ and $d\alpha_{F_2}$ are the same. It suffices then to show that
for any sufficiently small hypersurface $\Sigma \subset STM$ transverse to this common kernel line field, the restrictions of $d\alpha_{F_1}$ and $d\alpha_{F_2}$ coincide. 

By hypothesis, if $v_x$ is any given point in $STM$, the Reeb orbits for $\alpha_{F_1}$ and $\alpha_{F_2}$ with initial condition $v_x$ both pass through an open set $U \subset STM$ where $d\alpha_{F_1}$ and $d\alpha_{F_2}$ are  equal. Let $\Sigma \subset STM$ be a transverse hypersurface passing through $v_x$ such that its image $\Sigma_1$ under the time-$t_1$ map of the Reeb flow of $\alpha_{F_1}$ lies entirely in $U$. The surfaces $\Sigma$ and $\Sigma_1$ fall under the hypothesis of Lemma~\ref{surfaces-of-section} for both $(STM,\alpha_{F_1})$ and $(STM,\alpha_{F_2})$. Moreover, since  the map $\tau: \Sigma \rightarrow \Sigma_1$ described in Lemma~\ref{surfaces-of-section} 
depends only on the foliation defined by the Reeb orbits and not on their parametrization, it is
the same in both cases. It follows that $\tau$ is both an isomorphism between $(\Sigma, d\alpha_{F_1})$  and $(\Sigma_1, d\alpha_{F_1})$, and between $(\Sigma, d\alpha_{F_2})$ and $(\Sigma_1, d\alpha_{F_2})$. Since $d\alpha_{F_1}$ and $d\alpha_{F_2}$ are equal on  $\Sigma_1 \subset U$, we have that $d\alpha_{F_1}$ and $d\alpha_{F_2}$ are identical on $\Sigma$. 
\qed 

As an immediate consequence we recover a result from~\cite{Alvarez-Barbosa:2018} stating that solutions of Hilbert's fourth problem on open convex subsets  of $\R^n$ are determined, modulo a differential, by what happens ``at infinity". 

\begin{corollary}
If two $C^2$ projective Finsler metrics defined on an open convex subset $\mathcal{O} \subset \mathbb{R}^n$ agree outside of some compact set, then their difference is the differential of a $C^3$ function on $\mathcal{O}$. 
\end{corollary}

\medskip
We end the section with a modern formulation, due to Tabachnikov~\cite{Tabachnikov:1999a}, of a key
result in geodesic rigidity which goes back at least to Darboux~\cite[pp.~53--54, Livre~VI, Chapitre~III ]{Darboux:1894}. Since the standard proof for this result involves taking a Lie derivative of the differential of the Hilbert $1$-form---and thus requires that the metric be 
$C^3$---, we modify it slightly using Lemma~\ref{surfaces-of-section} so that it works for $C^2$ metrics. 

\begin{proposition}\label{integral-of-motion}
If $F_1$ and $F_2$ are two $C^2$ Finsler metrics on an $n$-dimen\-sional manifold $M$ and they have the same unparameterized geodesics, then there exists a  continuous, nowhere-zero function $\nu : \slit  \rightarrow \R$ that is homogeneous of degree zero and such that  $d\alpha_{{}_{F_2}}^{n-1} = \nu \, d\alpha_{{}_{F_1}}^{n-1}$. Moreover, $\nu$ is constant along the (common) geodesics of  $(M,F_1)$ and $(M,F_2)$. 
\end{proposition}

\proof
Since the forms $d\alpha_{{}_{F_1}}^{n-1}$ and $ d\alpha_{{}_{F_2}}^{n-1}$ are nowhere-zero $(2n-2)$-forms with the same kernel on the $(2n-1)$-dimensional space of tangent directions $STM = (\slit )/\R^+$, there exists a continuous, nowhere-zero function $\nu$ on $STM$ such that $d\alpha_{{}_{F_2}}^{n-1} = \nu d\alpha_{{}_{F_1}}^{n-1}$. Extending $\nu$ to $\slit$ as a homogeneous function of degree zero proves the first part of the proposition.

In order to show that $\nu$ is constant along orbits of the geodesic flow, take any  hypersurface $\Sigma \subset STM$ transverse to the Reeb flow of $\alpha_{{}_{F_1}}$ and let $\Sigma_t$ be
its image under the time-$t$ map of that flow. Since the Reeb foliations for $\alpha_{{}_{F_1}}$
and $\alpha_{{}_{F_2}}$ are the same, the map $\tau : \Sigma \rightarrow \Sigma_t$ defined in
Lemma~\ref{surfaces-of-section} is the same for both contact forms. It follows that 
$$
\tau^* d\alpha_{{}_{F_1}}^{n-1} = d\alpha_{{}_{F_1}}^{n-1} \textrm{ and that } 
\tau^* d\alpha_{{}_{F_2}}^{n-1} = d\alpha_{{}_{F_2}}^{n-1} .
$$
This implies that $\nu \circ \tau = \nu$. Since the transverse hypersurface $\Sigma$ and the time
$t$ are arbitrary, this shows that $\nu$ is constant on Reeb orbits. 
\qed

%%%%%%%%%%%%%%%%%%%%%%%%%%%%%%%%%%%%%%%%%%%%%%%%%%%%%%%%%%%%%%%%%%%
\section{Geodesically reversible Zoll metrics} \label{Zoll-metrics}
%%%%%%%%%%%%%%%%%%%%%%%%%%%%%%%%%%%%%%%%%%%%%%%%%%%%%%%%%%%%%%%%%%%

\begin{definition}
A Finsler manifold is {\sl Zoll} if its geodesic flow is periodic and its prime geodesics all have
the same length. Equivalently, a Finsler manifold is Zoll if its unit co-sphere bundle is a {\it regular contact manifold} in the sense of Boothby and Wang~\cite{Boothby-Wang:1958}.
\end{definition}

In Riemannian geometry, Zoll manifolds are very rare. By a theorem of Berger and Kazdan (see~\cite[pp.~236--246, Appendices~D and~E]{Besse:1978}) the only Zoll Riemannian metric in
$\RP^n$ is the standard constant curvature metric. It is not known whether there are Zoll Riemannian
metrics in $\CP^n$ or $\HP^n$ for $n > 1$ besides the standard metrics, and while there are many examples of Zoll metrics on spheres (see~\cite{Besse:1978}), they are not very explicit nor simple to describe. On the other hand there is an abundance of Zoll Finsler metrics. Non-reversible examples are in fact very cheap: take the standard metric in $S^n$, $\RP^n$, or
$\CP^n$ and compose the Hamiltonian of the metric with a homogeneous symplectic transformation that is sufficiently $C^2$ close to the identity. The result will be a (Zoll) Finsler metric whose geodesic flow is symplectomorphic to that of the standard metric.  Reversible examples of Zoll metric take a bit more work, but as we shall see they are still plentiful. 

\medskip
\noindent
{\it Reversible projective Finsler metrics on $S^n$ and $\RP^n$.} In his work on Hilbert's fourth
problem, Busemann had the lovely idea of inverting the classical Crofton formula for the length of curves on $\R^n$ and projective space: consider a smooth, positive measure $\mu$ on the space of projective hyperplanes in $\RP^n$, which is the dual projective space $\RP^{n*}$, and define the length of a $C^1$ curve by
the formula
$$
\length(\gamma) := \int_{H \in \RP^{n*}} \#(H \cap \gamma) \, d\mu .
$$
By construction, projective line segments are shortest curves and it can be shown that this definition of length corresponds to a smooth Finsler metric (see~\cite{Szabo:1986} and~\cite{Alvarez-Fernandes:2007}). These metrics can be lifted to the sphere $S^n$ and, indeed, every Finsler metric on $S^n$ for which great circles are geodesics is such a lift. Note that the space of metrics defined by this construction is parametrized by a smooth function---the density for the measure $\mu$---of $n$ variables. 

Note that if the measure $\mu$ is taken to be equal to the standard $SO(n+1)$-invariant measure everywhere except in the neighborhood of a point $H$ in $\RP^{n*}$, we obtain a Finsler metric that agrees with the standard metric outside a neighborhood of the projective hyperplane in $\RP^n$ defined by $H$. This provides us with examples of distinct Finsler metrics that have the same geodesics and agree on an open set. This is 
a Finsler phenomenon that has no Riemannian counterpart (see~\cite{Matveev-Topalov:2003}).

\medskip
\noindent
{\it Circular metrics on $\CP^n$ and $\HP^n$.} In~\cite{Alvarez-Duran:2002} \'Alvarez Paiva and 
Dur\'an Fernandez introduced Finsler submersions and used them to construct examples of metrics
on $\CP^n$ and $\HP^n$ where complex and quaternionic projective lines are totally geodesic and all
geodesics are (geometric) circles. This was done by using the idea of Busemann sketched above
to construct projective metrics on $S^{2n+1}$ and $S^{4n + 3}$ that are invariant under the Hopf actions of $S^1$ and $S^3$ on these spaces and then projecting those metrics down to
$\CP^n$ and $\HP^n$ by the Hopf map. Note that in the case of $\CP^n$ the space of these metrics
is parametrized by a  smooth function of $2n$ variables (i.e., $2n+1$ for a projective metric on $S^{2n +1}$ minus one variable because of $S^1$ invariance), and in the case of $\HP^n$ by a smooth function of $4n$ variables. 

\medskip
\noindent
{\it Circular metrics on $S^2$.} In~\cite{Alvarez-Berk:2010} \'Alvarez Paiva and Berck constructed 
all reversible Finsler metrics on the two-sphere for which geodesics are geometric circles. The construction depends on two functions of two variables: one function to define the circular path geometry and another to metrize it. 

\medskip
\noindent
{\it Smooth projective planes.} A smooth real projective plane is a reversible path geometry on $\RP^2$ such that any two distinct points determine a unique path and any two paths intersect in only one point. Moreover, the path joining two distinct points and the intersection point of two distinct paths are required to depend smoothly on the points and paths, respectively. 

In order to obtain a great number of examples of these planes consider the point-line incidence relation 
$$
I := \{(p,l) \in \RP^2 \times \RP^{2*} : \textrm{ the point $p$ lies on the line $l$ } \} .
$$
As remarked by McKay (Corollary~1 in~\cite{McKay:2005}), it follows from the work of B\"odi and Immervoll~(\cite{Bodi-Immervoll:2000}) that any sufficently small $C^1$ perturbation $\Tilde{I}$ of the point-line incidence relation as a submanifold of $\RP^2 \times \RP^{2*}$ defines a smooth projective plane:  the path in $\RP^2$ corresponding to a point $l \in \RP^{2*}$ is the 
set of points $p \in \RP^2$ such that $(p,l) \in \Tilde{I}$.

There are even examples of smooth projective planes, due to Segre and Immervoll (see~\cite{Immervoll:2003}) that are real analytic. By Theorem~2.2 in~\cite{Alvarez-Berk:2010} every smooth projective plane can be metrized  as smooth reversible Finsler manifolds. 

After this short review of all known constructions of non-Riemannian reversible Zoll Finsler metrics, we are now ready to pass to the main results of this section.  

\begin{maintheorem}{II}
If $(M,F)$ is a smooth, geodesically-reversible Zoll Finsler manifold, then $F$ is the sum of a reversible Zoll Finsler metric and the differential of a smooth function on $M$.
\end{maintheorem}

Before going into the proof of this result, we need a couple of simple lemmas. The first is 
a direct consequence of the results of Boothby and Wang~\cite{Boothby-Wang:1958} on regular contact manifolds (see~\cite{Weinstein:1974}  or Proposition~3.3 in \cite{Alvarez-Balacheff:2014}).

\begin{lemma}\label{Zoll-volume}
The volume of an $n$-dimensional Zoll Finsler manifold is equal to an integer multiple of 
$(\varepsilon_n n!)^{-1}$ times the $n$-th power of the length of (any) one of its prime geodesics.
\end{lemma}

The second lemma was worked out in a MathOverflow exchange (see~\cite{Alvarez-MathOverflow:2013}). 

\begin{lemma}\label{Fundamental-group-Zoll}
The order of the fundamental group of a geodesically reversible Zoll Finsler manifold $(M,F)$  is at most two. In particular, every closed differential $1$-form on $M$ is exact. 
\end{lemma}

\proof
 Assume $M$ is not simply connected and let $\sigma \subset M$ be a non-null-homotopic closed curve. By a classic theorem of Hilbert $\sigma$ is homotopic to a closed geodesic and that geodesic is covered a certain number of times $n \in \Z$ by a prime geodesic $\gamma$. Remark now that any two prime geodesics $\gamma_1$ and $\gamma_2$ on a Zoll manifold are homotopic: if $v_1$ is a velocity vector of $\gamma_1$ and $v_2$ is a velocity vector of $\gamma_2$ any continuous path $v_t$ in $\slit$ from $v_1$ to $v_2$ defines a homotopy between $\gamma_1$ and $\gamma_2$ by associating to $v_t$ the unique prime geodesic tangent to $v_t$. This shows that the fundamental group of $M$ has only one generator. Moreover, since $(M,F)$ is geodesically reversible, the 
 reverse of any prime geodesic $\gamma$ is also a prime geodesic, which means that $\gamma$ is homotopic to its reverse and hence the fundamental group of $M$ has order two.
\qed

\proof[Proof of Theorem~\ref{Zoll}]
Note that since $F$ is geodesically reversible, the reverse metric $F_{-}$ has the same unparametrized geodesics. Moreover, its prime geodesics have the same length as those of $F$.
It follows that the convex combination  $(1-t)F + tF_{-}$ is a continuous path of Zoll metrics
whose prime geodesics have the same length. By Lemma~\ref{Zoll-volume} the volume of all
these metrics must be the same.  However, note that for $t = 1/2$ we have the central symmetrization
$F^S$ whose volume, by Proposition~\ref{central-symmetrization}, must be greater than or equal to that of $F$ with equality only if $F = F_S + \beta$, where $\beta$ is a $1$-form. By Corollary~\ref{adding-form}, $\beta$ must be closed and Lemma~\ref{Fundamental-group-Zoll} implies it must be exact. 
\qed

Using Theorem~\ref{Zoll} we can now settle Hilbert's fourth problem for asymmetric metrics on the $n$-sphere. 

\begin{maintheorem}{I}
 If $d$ is a continuous asymmetric metric on the $n$-sphere $S^n$ for which all great circles are geodesics, then there exists a (continuous) function $f : S^n \rightarrow \R$ such that
$d(x,y) - d(y,x) = f(y) - f(x)$.
\end{maintheorem}

\proof
By a result of Pogorelov~\cite{Pogorelov:1978} (see also~\cite{Szabo:1986}) every continuous
metric $d$ on the $n$-sphere for which all great circles are geodesics is an uniform limit of distance functions $d_n$ induced by smooth projective Finsler metrics. In the cited papers the 
result is stated for reversible metrics, but this hypothesis is never used in the proof. 

By Theorem~\ref{Zoll}, for each $d_n$ we have $d_n(x,y) - d_n(y,x) = f_n(y) - f_n(x)$ for some smooth function $f_n$ on the $n$-sphere. Since the functions $f_n$ are defined up to an additive constant, we may fix an arbitrary point $\zeta$ on $S^n$ and assume that all the $f_n$ vanish at that point. Define $f$ as the limit
$$
f(x) :=  \lim_{n \rightarrow \infty} \left( d_n(\zeta,x) - d_n(x,\zeta) \right) ,
$$
which, on account of the uniform convergence of the $d_n$, is a continuous function. It follows that 
\begin{align*}
 d(x,y) - d(y,x) &= \lim_{n \rightarrow \infty} \left( d_n(x,y) - d_n(y,x) \right) = 
 \lim_{n \rightarrow \infty} \left( f_n(y) - f_n(x)\right) \\ 
 &= f(y) - f(x) .   
\end{align*}
\qed

The ideas in the proof of Theorem~\ref{Zoll} unveil a second interesting feature of solutions to Hilbert's fourth problem on $\RP^n$: up to a differential, they are determined by their Holmes-Thompson volume densities. In order to prove this, we will need a sharper version of Lemma~\ref{Zoll-volume} due to Weinstein~\cite{Weinstein:1978} and  Yang~\cite{Yang:1980} and which was a key component of the proof of the Blaschke conjecture. 

\begin{theorem}[Weinstein and Yang]\label{weak-Blaschke-conjeture}
If $F$ is a $C^2$ Zoll Finsler metric on $S^n$ whose prime geodesics have length 
$2\pi r$, then its Holmes-Thompson volume is the same as the volume of the round sphere of radius $r$. 
\end{theorem}

In~\cite{Weinstein:1974} and~\cite{Yang:1980} Weinstein and Yang state their theorems for Riemannian metrics, however the proof, which just relies on the algebraic topology of the fibrations involving the unit sphere bundle, the $n$-sphere, and the manifold of geodesics, holds without modification for Finsler metrics.  In what follows we will need only the following consequence of their results:

\begin{lemma}\label{Projective-metric-volume}
Two $C^2$ projective Finsler metrics on $\RP^n$ have the same Holmes-Thompson volume if and only if their prime geodesics have the same length. 
\end{lemma}

\proof
Lift the metrics to $S^n$, the double cover of $\RP^n$, and apply Theorem~\ref{weak-Blaschke-conjeture}.
\qed

\begin{theorem}
The difference of two $C^2$ projective Finsler metrics on $\RP^n$ ($n > 1$) that have the same Holmes-Thompson density is the differential of a $C^3$ function.
\end{theorem}

\proof[Sketch of the proof.]
If two projective metrics $F_1$ and $F_2$ on $\RP^n$ have the same Holmes-Thompson volume (let alone identical volume densities), then Lemma~\ref{Projective-metric-volume} tells us that their prime geodesics have the same length, $\ell$. The average metric $(F_1 + F_2)/2$ is a projective Fisler metric with prime geodesics of length $\ell$ and, therefore, its Holmes-Thompson volume equal to those of $F_1$ and $F_2$. However, by Proposition~\ref{Brunn-Minkowski}, this can only be the case if  $F_2 = F_1 + \beta$ for some $1$-form $\beta$. As in the proof of Theorem~\ref{Zoll}, we conclude that $\beta$ is exact. 
\qed
 
%%%%%%%%%%%%%%%%%%%%%%%%%%%%%%%%%%%%%%%%%%%%%%%%%%%%%%%%%%%%%%%%%%%%%%%%%%%%%%%%%%%%%% 

%%%%%%%%%%%%%%%%%%%%%%%%%%%%%%%%%%%%%%%%%%%%%%%%%%%%%%%%%%%%%%%%%%%%%%%%%%%%%%%%%%%%%%%%%%%%
\section{Finsler metrics without non-trivial integrals of motion} \label{geodesic-rigidity}
%%%%%%%%%%%%%%%%%%%%%%%%%%%%%%%%%%%%%%%%%%%%%%%%%%%%%%%%%%%%%%%%%%%%%%%%%%%%%%%%%%%%%%%%%%%%
Every continuous function on the space of geodesics of a Zoll manifold $(M,F)$ is an integral of motion for the geodesic flow of the metric $F$. At the other end of the spectrum we have Finsler manifolds where the geodesic flow is Anosov or transitive and the only integrals of motion are trivial: constant multiples of the Hamiltonian of the metric. The reader will find many classes of geometrically meaningful examples in the works of Foulon~\cite{Foulon:1986}, Barbosa Gomes, Ruggeiro and Chimenton~\cite{Barbosa-Ruggiero:2013, Barbosa-Chimenton-Ruggiero:2020}. 

In two dimensions, geodesic rigidity theorems for metrics without non-trivial integrals of motion are very easy to prove:  if $F_1$ and $F_2$ have the same geodesics and no non-trivial integrals of motion, then by Proposition~\ref{integral-of-motion} the differentials of the Hilbert $1$-forms are
multiples and therefore one metric is a multiple of the other plus a closed $1$-form. The main problem in $n$ dimensions ($n > 2$) is to extract meaningful geometric information from the proportionality of $d\alpha_{F_1}^{n-1}$ and $d\alpha_{F_2}^{n-1}$. This difficulty is solved 
by a simple application of the theory of areas in Finsler spaces. 

\begin{theorem} \label{areas}
Let $F_1$ and $F_2$ be two $C^2$ Finsler metrics defined on a manifold $M$ which have the same unparametrized geodesics. If the geodesic flow of either metric does not admit non-trivial 
integrals of motion, there is a constant $\nu > 0$ such that the Holmes-Thompson area of any hypersurface in $(M,F_2)$ is $\nu$ times its area in  $(M,F_1)$.
\end{theorem}

\proof
Since $F_1$ and $F_2$ have the same unparametrized geodesics, Proposition~\ref{integral-of-motion}
tells us that there exists a positive, continuous function $\nu$ on the manifold of tangent directions $STM$ that is constant along geodesics and such that  $d\alpha_{F_2}^{n-1} = \nu d\alpha_{F_1}^{n-1}$. By hypothesis,  $\nu$ is a constant. 

 Let $N \subset M$ be a hypersurface, which, by the locality of area, we can assume to be contained in a region $B \subset M$ that is convex for both $F_1$ and $F_2$. The set of (common) geodesics in $B$ is a manifold $\geod{B}$ and we denote by $\omega_1$ and $\omega_2$ the maximally-non-degenerate, continuous $2$-forms induced by the differentials of the Hilbert forms $\alpha_{F_1}$ and $\alpha_{F_2}$ in $STM$. The equality $d\alpha_{F_2}^{n-1} = \nu d\alpha_{F_1}^{n-1}$ translates to the equality  $\omega_2^{n-1} = \nu \omega_1^{n-1}$, which in turn yields
$$
 \int_{\gamma \in \geod{K}}
\#(\gamma \cap N) |\omega_2^{n-1}| = 
\nu \int_{\gamma \in \geod{K}}
\#(\gamma \cap N) |\omega_1^{n-1}| .
$$ 
By the Crofton formula for hypersurfaces in Finsler spaces (Theorem~\ref{crofton-hypersurface}) the area of $N$ in $(M,F_2)$ equals $\nu$ times the area of $N$ in $(M,F_1)$.
\qed

\begin{maintheorem}{VI} 
Let $F_1$ and $F_2$ be two reversible $C^2$ Finsler metrics defined on a manifold $M$ and with the same unparametrized geodesics. If the geodesic flow does not admit any non-trivial integral of motion, then the metric $F_2$ is necessarily a  constant multiple of $F_1$.
\end{maintheorem}

\proof
Reversible Finsler metrics are completely determined by their hypersurface area density (Proposition~\ref{areas-determine-metric}), therefore Theorem~\ref{areas} implies that $F_2$ is a constant multiple of $F_1$. 
\qed

\begin{maintheorem}{V}
Let $M$ be a two- or three-dimensional manifold and let $F$ be a geodesically reversible Finsler metric on $M$. If the geodesic flow does not admit any non-trivial integral of motion, then the metric $F$ is necessarily the sum of a reversible metric and a closed 1-form.
\end{maintheorem}

\proof
Since the Finsler metric $F$ is geodesically reversible it has the same un\-parametrized geodesics as its central symmetrization $F^S$. Theorem~\ref{areas} tells us that there exists a constant $\nu > 0$ such that the hypersurface area density of $(M,F^S)$ equals $\nu$ times that of $(M,F)$. By definition, the hypersurface area density of the metric $F$ equals that of its areal symmetrization $F^B$ and, since both $F^S$ and $F^B$ are symmetric, Proposition~\ref{areas-determine-metric} implies that  $F^S$ and $F^B$ are proportional and, a fortiori, conformal.  Proposition~\ref{constant-width-constant-brightness} now implies that $F = F_S + \beta$, where $\beta$ is a $1$-form, which is closed by Proposition~\ref{adding-form}.
\qed

Theorems~\ref{hyperbolic-space-forms-non-reversible} and~\ref{hyperbolic-space-forms-reversible}
follow immediately from the previous results by recalling that the geodesic flow on a hyperbolic space form of finite volume is transitive.

%%%%%%%%%%%%%%%%%%%%%%%%%%%%%%%%%%%%%%%%%%%%%%%%%%%%%%%%%%%%%%%%%%%%%%
%%%%%%%%%%%%%%%%%%%%%%%%%%%%%% REFERENCES %%%%%%%%%%%%%%%%%%%%%%%%%%%%
%%%%%%%%%%%%%%%%%%%%%%%%%%%%%%%%%%%%%%%%%%%%%%%%%%%%%%%%%%%%%%%%%%%%%%

\bibliography{paperbib}
\bibliographystyle{amsplain}

%%%%%%%%%%%%%%%
\end{document}